%&amstex
\input amstex
\input amsppt.sty
\magnification=\magstep1
\hsize=30truecc
\baselineskip=16truept
\vsize=22.2truecm
%\NoBlackBoxes
\nologo
\pageno=1
\topmatter
\TagsOnRight

\def\Z{\Bbb Z}

\def\N{\Bbb N}

\def\al{\alpha}
\def\l{\left}
\def\r{\right}
\def\bg{\bigg}
\def\({\bg(}
\def\){\bg)}
\def\[{\bg[}
\def\]{\bg]}
\def\t{\text}
\def\f{\frac}

\def\eq{\equiv}

\def\ls{\leqslant}
\def\gs{\geqslant}
\def\mo{\t{\rm mod}}

\def\Proof{\noindent{\it Proof}}

\def\Remark{\medskip\noindent{\it  Remark}}
\def\Ack{\medskip\noindent {\bf Acknowledgment}}
\hbox{Acta Arith. 127(2007), no.\,2, 103--113.}
\bigskip
\title Mixed sums of squares and triangular numbers
\endtitle
\author {Zhi-Wei Sun (Nanjing)}\endauthor
\leftheadtext{Z. W. Sun}

\abstract  By means of $q$-series, we prove that
any natural number
is a sum of an even square and two triangular numbers,
and that each positive integer is a sum of a triangular number
plus $x^2+y^2$ for some $x,y\in\Z$ with
$x\not\eq y\ (\mo\ 2)$ or $x=y>0$.
The paper also contains some other results and open conjectures
on mixed sums of squares and triangular numbers.
\endabstract

\thanks 2000 {\it Mathematics Subject Classification}:\,Primary 11E25;
Secondary 05A30, 11D85, 11P99.
\newline\indent Supported by the National Science Fund for
Distinguished Young Scholars (Grant No. 10425103) in China.
\endthanks
\endtopmatter
\document

\heading{1. Introduction}\endheading

 A classical result of Fermat asserts that any prime $p\eq 1\ (\mo\ 4)$
 is a sum of two squares of integers. Fermat also conjectured
 that each $n\in\N$ can be written as a sum of three triangular numbers, where $\N$
 is the set $\{0,1,2,\ldots\}$ of natural numbers,
 and triangular numbers are those integers $t_x=x(x+1)/2$ with $x\in\Z$.
 An equivalent version of this conjecture states that
 $8n+3$ is a sum of three squares (of odd integers).
 This follows from the following profound theorem (see, e.g.,
 [G, pp.\,38--49] or [N, pp.\,17-23]).

 \proclaim{Gauss-Legendre Theorem} $n\in\N$ can be written as a sum of three squares
 of integers
 if and only if $n$ is not of the form $4^k(8l+7)$ with $k,l\in\N$.
 \endproclaim

 Building on some work of Euler, in 1772 Lagrange showed that
 every natural number is a sum of four squares of integers.

 For problems and results on representations of natural numbers by various
 quadratic forms with coefficients in $\N$, the reader may consult [Du] and [G].

 Motivated by Ramanujan's work [Ra],
 L. Panaitopol [P] proved the following interesting result in 2005.

 \proclaim{Theorem A} Let $a,b,c$ be positive integers
 with $a\ls b\ls c$.
 Then every odd natural number can be written in the form $ax^2+by^2+cz^2$ with $x,y,z\in\Z$,
 if and only if the vector $(a,b,c)$ is $(1,1,2)$ or $(1,2,3)$ or $(1,2,4)$.
 \endproclaim

 According to L. E. Dickson [D2, p.\,260], Euler already noted that
 any odd integer $n>0$ is representable by $x^2+y^2+2z^2$ with $x,y,z\in\Z$.

 In 1862 J. Liouville (cf. [D2, p.\,23]) proved the following result.

 \proclaim{Theorem B} Let $a,b,c$ be positive integers with $a\ls b\ls c$.
Then every $n\in\N$ can be written as $at_x+bt_y+ct_z$ with
$x,y,z\in\Z$, if and only if $(a,b,c)$ is among the following
vectors:
$$(1,1,1),\ (1,1,2),\ \ (1,1,4),\ (1,1,5),\ (1,2,2),
\ (1,2,3),\ (1,2,4).$$
\endproclaim

 Now we turn to representations of natural numbers
 by mixed sums of squares (of integers)
 and triangular numbers.

 Let $n\in\N$. By the Gauss-Legendre theorem, $8n+1$ is a sum of three squares.
 It follows that $8n+1=(2x)^2+(2y)^2+(2z+1)^2$
for some $x,y,z\in\Z$ with $x\eq y\ (\mo\ 2)$;
this yields the representation
$$n=\f{x^2+y^2}2+t_z=\l(\f{x+y}2\r)^2+\l(\f{x-y}2\r)^2+t_z$$
as observed by Euler. According to Dickson [D2, p.\,24], E.
Lionnet stated, and V. A. Lebesgue [L] and M. S. R\'ealis [Re]
proved that $n$ can also be written in the form $x^2+t_y+t_z$ with
$x,y,z\in\Z$. Quite recently, this was reproved by H. M. Farkas
[F] via the theory of theta functions.

 Using the theory of ternary quadratic forms,
in 1939 B. W. Jones and G. Pall [JP, Theorem 6]
proved that for any $n\in\N$ we have $8n+1=ax^2+by^2+cz^2$
for some $x,y,z\in\Z$ if the vector $(a,b,c)$ belongs to the set
$$\{(1,1,16),\ (1,4,16),\ (1,16,16),\ (1,2,32),\ (1,8,32),\ (1,8,64)\}.$$
As $(2z+1)^2=8t_z+1$, the result of Jones and Pall implies that
each $n\in\N$ can be written in any of the following three forms
with $x,y,z\in\Z$:
$$2x^2+2y^2+t_z=(x+y)^2+(x-y)^2+t_z,\ x^2+4y^2+t_z,\ x^2+8y^2+t_z.$$

In this paper we establish the following result by means of $q$-series.

\proclaim{Theorem 1} {\rm (i)} Any $n\in\N$ is a sum of an even square and two triangular numbers.
Moreover, if $n/2$ is not a triangular number then
$$\aligned&|\{(x,y,z)\in\Z\times\N\times\N:\ x^2+t_y+t_z=n\ \t{and}\ 2\nmid x\}|
\\=&|\{(x,y,z)\in\Z\times\N\times\N:\ x^2+t_y+t_z=n\ \t{and}\ 2\mid x\}|.
\endaligned\tag1$$

{\rm (ii)} If $n\in\N$ is not a triangular number, then
$$\aligned&|\{(x,y,z)\in\Z\times\Z\times\N:\ x^2+y^2+t_z=n\ \t{and}\ x\not\eq y\ (\mo\ 2)\}|
\\=&|\{(x,y,z)\in\Z\times\Z\times\N:\ x^2+y^2+t_z=n\ \t{and}\ x\eq y\ (\mo\ 2)\}|>0.
\endaligned\tag2$$

{\rm (iii)} A positive integer $n$ is a sum of an odd square, an
even square and a triangular number, unless it is a triangular
number $t_m\ (m>0)$ for which all prime divisors of $2m+1$ are
congruent to $1$ mod $4$ and hence $t_m=x^2+x^2+t_z$ for some
integers $x>0$ and $z$ with $x\eq m/2\ (\mo\ 2)$.
\endproclaim

\Remark. Note that $t_2=1^2+1^2+t_1$ but we cannot write $t_2=3$
as a sum of an odd square, an even square and a triangular number.
\medskip

Here are two more theorems of this paper.

\proclaim{Theorem 2} Let $a,b,c$ be positive integers with $a\ls
b$. Suppose that every $n\in\N$ can be written as $ax^2+by^2+ct_z$
with $x,y,z\in\Z$. Then $(a,b,c)$ is among the following vectors:
$$\gather(1,1,1),\ (1,1,2),\ (1,2,1),\ (1,2,2),\ (1,2,4),
\\(1,3,1),\ (1,4,1),\ (1,4,2),\ (1,8,1),\ (2,2,1).
\endgather$$
\endproclaim

\proclaim{Theorem 3} Let $a,b,c$ be positive integers with $b\gs
c$. Suppose that every $n\in\N$ can be written as $ax^2+bt_y+ct_z$
with $x,y,z\in\Z$. Then $(a,b,c)$ is among the following vectors:
$$\gather(1,1,1),\ (1,2,1),\ (1,2,2),\ (1,3,1),\ (1,4,1),\ (1,4,2),\ (1,5,2),
\\ (1,6,1),\ (1,8,1),
\ (2,1,1),\ (2,2,1),\ (2,4,1),\ (3,2,1),
\ (4,1,1),\ (4,2,1).
\endgather$$
\endproclaim

Theorem 1 and Theorems 2--3 will be proved in Sections 2 and 3 respectively.
In Section 4, we will pose three conjectures and discuss the converses of Theorems 2 and 3.

\heading{2. Proof of Theorem 1}\endheading

Given two integer-valued quadratic polynomials $f(x,y,z)$ and
$g(x,y,z)$, by $f(x,y,z)\sim g(x,y,z)$ we mean
$$\{f(x,y,z):\,x,y,z\in\Z\}=\{g(x,y,z):\,x,y,z\in\Z\}.$$
Clearly $\sim$ is an equivalence relation on the set of all
integer-valued ternary quadratic polynomials.

The following lemma is a refinement of Euler's observation
$t_y+t_z\sim y^2+2t_z$ (cf. [D2, p.\,11]).

\proclaim{Lemma 1} For any $n\in\N$ we have
$$|\{(y,z)\in\N^2:\, t_y+t_z=n\}|=|\{(y,z)\in\Z\times\N:\, y^2+2t_z=n\}|.\tag3$$
\endproclaim
\Proof. Note that $t_{-y-1}=t_y$. Thus
$$\align&|\{(y,z)\in\N^2:\, t_y+t_z=n\}|
=\f14|\{(y,z)\in\Z^2:\, t_y+t_z=n\}|
\\=&\f14|\{(y,z)\in\Z^2:\, 4n+1=(y+z+1)^2+(y-z)^2\}|
\\=&\f14|\{(x_1,x_2)\in\Z^2:\, 4n+1=x_1^2+x_2^2\ \t{and}\ x_1\not\eq x_2\ (\mo\ 2)\}|
\\=&\f24|\{(y,z)\in\Z^2:\, 4n+1=(2y)^2+(2z+1)^2\}|
\\=&\f12|\{(y,z)\in\Z^2:\, n=y^2+2t_z\}|
=|\{(y,z)\in\Z\times\N:\, n=y^2+2t_z\}|.
\endalign$$
This concludes the proof. \qed
\medskip

 Lemma 1 is actually equivalent to the following observation of Ramanujan (cf. Entry 25(iv) of [B, p.\,40]):
$\psi(q)^2=\varphi(q)\psi(q^2)$ for $|q|<1$, where
$$\varphi(q)=\sum_{n=-\infty}^{\infty}q^{n^2}\ \ \t{and}\ \
\psi(q)=\sum_{n=0}^{\infty}q^{t_n}.\tag4$$

Let $n\in\N$ and define
$$\gather r(n)=|\{(x,y,z)\in\Z\times\N\times\N:\ x^2+t_y+t_z=n\}|,\tag5
\\r_0(n)=|\{(x,y,z)\in\Z\times\N\times\N:\ x^2+t_y+t_z=n\ \t{and}\ 2\mid x\}|,\tag6
\\r_1(n)=|\{(x,y,z)\in\Z\times\N\times\N:\ x^2+t_y+t_z=n\ \t{and}\ 2\nmid x\}|.\tag7
\endgather$$
Clearly $r_0(n)+r_1(n)=r(n)$. In the following lemma we investigate the difference
$r_0(n)-r_1(n)$.

\proclaim{Lemma 2} For $m=0,1,2,\ldots$ we have
$$r_0(2t_m)-r_1(2t_m)=(-1)^m(2m+1).\tag 8$$
Also, $r_0(n)=r_1(n)$ if $n\in\N$ is not a triangular number times $2$.
\endproclaim
\Proof. Let $|q|<1$. Recall the following three known identities
implied by Jacobi's triple product identity
(cf. [AAR, pp.\,496--501]):
$$\align\varphi(-q)=&\prod_{n=1}^\infty(1-q^{2n-1})^2(1-q^{2n})\quad\t{(Gauss)},
\\\psi(q)=&\prod_{n=1}^\infty\f{1-q^{2n}}{1-q^{2n-1}}\quad\t{(Gauss)},
\\\prod_{n=1}^\infty(1-q^n)^3=&\sum_{n=0}^\infty(-1)^n(2n+1)q^{t_n}\quad\t{(Jacobi)}.
\endalign$$
Observe that
$$\align&\sum_{n=0}^{\infty}(r_0(n)-r_1(n))q^n
\\=&\(\sum_{x=-\infty}^{\infty}(-1)^xq^{x^2}\)
\(\sum_{y=0}^{\infty}q^{t_y}\)\(\sum_{z=0}^{\infty}q^{t_z}\)=\varphi(-q)\psi(q)^2
\\=&\(\prod_{n=1}^{\infty}(1-q^{2n-1})^2(1-q^{2n})\)
\(\prod_{n=1}^\infty\f{1-q^{2n}}{1-q^{2n-1}}\)^2
\\=&\prod_{n=1}^\infty(1-q^{2n})^3=\sum_{m=0}^\infty(-1)^m(2m+1)(q^2)^{t_m}.
\endalign$$
Comparing the coefficients of $q^n$ on both sides, we obtain the
desired result. \qed
\medskip

The following result was discovered by Hurwitz in 1907 (cf. [D2,
p.\,271]); an extension was established in [HS] via the theory of
modular forms of half integer weight.

\proclaim{Lemma 3} Let $n>0$ be an odd integer, and let $p_1,\ldots,p_r$ be all
the distinct prime divisors of $n$
congruent to $3$ mod $4$.
Write $n=n_0\prod_{0<i\ls r}p_i^{\al_i}$,  where $n_0,\al_1,\ldots,\al_r$ are positive
integers and $n_0$ has no prime divisors congruent to $3$ mod $4$.
Then
$$|\{(x,y,z)\in\Z^3:\ x^2+y^2+z^2=n^2\}|
=6n_0\prod_{0<i\ls r}\l(p_i^{\al_i}+2\f{p_i^{\al_i}-1}{p_i-1}\r).\tag9$$
\endproclaim
\Proof. We deduce (9) in a new way and use some standard notations in number theory.

 By (4.8) and (4.10) of [G],
 $$\align&|\{(x,y,z)\in\Z^3:\, x^2+y^2+z^2=n^2\}|
 \\=&\sum_{d\mid n}\f{24}{\pi}d\sum_{m=1}^\infty\f1m\l(\f{-4d^2}m\r)
 \\=&\f{24}{\pi}\sum_{d\mid n}d\sum_{k=1}^\infty\f{(-1)^{k-1}}{2k-1}\sum_{c\mid \gcd(2k-1,d)}\mu(c)
 \\=&\f{24}{\pi}\sum_{d\mid n}d\sum_{c\mid d}\mu(c)\sum_{k=1}^\infty\f{(-1)^{((2k-1)c-1)/2}}{(2k-1)c}
 \\=&\f{24}{\pi}\sum_{d\mid n}d\sum_{c\mid d}\f{\mu(c)}c(-1)^{(c-1)/2}\sum_{k=1}^\infty\f{(-1)^{k-1}}{2k-1}
 \\=&6\sum_{c\mid n}(-1)^{(c-1)/2}\f{\mu(c)}c\sum_{q\mid\f nc}cq
 =6\sum_{c\mid n}(-1)^{(c-1)/2}\mu(c)\sigma\l(\f nc\r)
 \endalign$$
and hence
$$\align&\f16|\{(x,y,z)\in\Z^3:\, x^2+y^2+z^2=n^2\}|
\\=&\sum_{d_0\mid n_0}\ \sum_{d_1\mid p_1^{\al_1},\ldots,d_r\mid p_r^{\al_r}}\l(\f{-1}{d_0d_1\cdots d_r}\r)
 \mu(d_0d_1\cdots d_r)\sigma\l(\f{n_0}{d_0}\prod_{0<i\ls r}\f{p_i^{\al_i}}{d_i}\r)
 \\=&\sum_{d_0\mid n_0}\l(\f{-1}{d_0}\r)\mu(d_0)\sigma\l(\f{n_0}{d_0}\r)
 \times\prod_{0<i\ls r}\sum_{d_i\mid p_i^{\al_i}}\l(\f{-1}{d_i}\r)\mu(d_i)\sigma\l(\f{p_i^{\al_i}}{d_i}\r)
\\=&\sum_{d_0\mid n_0}\mu(d_0)\sigma\l(\f{n_0}{d_0}\r)
 \times\prod_{0<i\ls r}\l(\sigma(p_i^{\al_i})+\l(\f{-1}{p_i}\r)\mu(p_i)\sigma(p_i^{\al_i-1})\r)
 \\=&n_0\prod_{0<i\ls r}\l(p_i^{\al_i}+2\sigma(p_i^{\al_i-1})\r)
 =n_0\prod_{0<i\ls r}\l(p_i^{\al_i}+2\f{p_i^{\al_i}-1}{p_i-1}\r).
\endalign$$
This completes the proof. \qed

\medskip
\noindent{\tt Proof of Theorem 1}. (i) By the Gauss-Legendre
theorem, $4n+1$ is a sum of three squares and hence
$4n+1=(2x)^2+(2y)^2+(2z+1)^2$ (i.e., $n=x^2+y^2+2t_z$) for some
$x,y,z\in\Z$. Combining this with Lemma 1 we obtain a simple proof
of the known result
 $$r(n)=|\{(x,y,z)\in\Z\times\N\times\N:\,x^2+t_y+t_z=n\}|>0.$$

Recall that $r_0(n)+r_1(n)=r(n)$. If $n/2$ is not a triangular
number, then $r_0(n)=r_1(n)=r(n)/2>0$ by Lemma 2. If $n=2t_m$ for
some $m\in\N$, then we also have $r_0(n)>0$ since $n=0^2+t_m+t_m$.
\smallskip

(ii) Note that
$$n=x^2+y^2+t_z \iff 2n=2(x^2+y^2)+2t_z=(x+y)^2+(x-y)^2+2t_z.$$
From this and Lemma 1, we get
$$\align&|\{(x,y,z)\in\Z\times\Z\times\N:\, x^2+y^2+t_z=n\}|
\\=&|\{(x,y,z)\in\Z\times\Z\times\N:\, x^2+y^2+2t_z=2n\}|
\\=&|\{(x,y,z)\in\Z\times\N\times\N:\, x^2+t_y+t_z=2n\}|=r(2n)>0;
\endalign$$
in the language of generating functions, it says that
$$\varphi(q)\psi(q)^2+\varphi(-q)\psi(-q)^2=2\varphi(q^2)^2\psi(q^2).$$
Similarly,
$$\align&|\{(x,y,z)\in\Z\times\Z\times\N:\, x^2+y^2+t_z=n\ \t{and}\ x\eq y\ (\mo\ 2)\}|
\\=&|\{(x,y,z)\in\Z\times\Z\times\N:\, x^2+y^2+2t_z=2n\ \t{and}\ 2\mid x\}|
\\=&|\{(x,y,z)\in\Z\times\N\times\N:\, x^2+t_y+t_z=2n\ \t{and}\ 2\mid x\}|
=r_0(2n)>0
\endalign$$
and
$$|\{(x,y,z)\in\Z\times\Z\times\N:\, x^2+y^2+t_z=n\ \t{and}\ x\not\eq y\ (\mo\ 2)\}|=r_1(2n).\tag10$$
If $n$ is not a triangular number, then
$r_0(2n)=r_1(2n)=r(2n)/2>0$ by Lemma 2, and hence (2) follows from the above.
\smallskip

(iii) By Theorem 1(ii), if $n$ is not a triangular number then
$n=x^2+y^2+t_z$ for some $x,y,z\in\Z$ with $2\mid x$ and $2\nmid
y$.

Now assume that $n=t_m\ (m>0)$ is not a sum of an odd square, an
even square and a triangular number. Then $r_1(2t_m)=0$ by (10).
In view of (ii) and (8),
$$\align&|\{(x,y,z)\in\Z\times\Z\times\N:\,x^2+y^2+t_z=t_m\}|
\\=&r_0(2t_m)+r_1(2t_m)=r_0(2t_m)-r_1(2t_m)=(-1)^m(2m+1).
\endalign$$
Therefore
$$\align &(-1)^m(2m+1)=\f12|\{(x,y,z)\in\Z^3:\ x^2+y^2+t_z=t_m\}|
\\=&\f12|\{(x,y,z)\in\Z^3:\ 2(2x)^2+2(2y)^2+(2z+1)^2=8t_m+1\}|
\\=&\f12|\{(x,y,z)\in\Z^3: (2x+2y)^2+(2x-2y)^2+(2z+1)^2=8t_m+1\}|
\\=&\f12|\{(x_1,y_1,z)\in\Z^3:\ 4(x_1^2+y_1^2)+(2z+1)^2=(2m+1)^2\}|
\\=&\f16|\{(x,y,z)\in\Z^3:\ x^2+y^2+z^2=(2m+1)^2\}|.
\endalign$$

Since
$$(-1)^m6(2m+1)=|\{(x,y,z)\in\Z^3:\ x^2+y^2+z^2=(2m+1)^2\}|\not>6(2m+1),$$
by Lemma 3 the odd number $2m+1$ cannot have a prime divisor congruent to 3 mod 4.
So all the prime divisors of $2m+1$
are congruent to 1 mod 4, and hence
$$|\{(x,y)\in\N^2:\ x>0\ \t{and}\ x^2+y^2=2m+1\}|=\sum_{d\mid 2m+1}1>1$$
by Proposition 17.6.1 of [IR, p.\,279]. Thus $2m+1$ is a sum of two squares of positive integers.
Choose positive integers $x$ and $y$ such that $x^2+y^2=2m+1$ with $2\mid x$ and $2\nmid y$.
Then
$$8t_m+1=(2m+1)^2=(x^2-y^2)^2+4x^2y^2=8t_{(x^2-y^2-1)/2}+1+16\l(\f x2\r)^2y^2$$
and hence $$t_m=t_{(x^2-y^2-1)/2}+\l(\f x2y\r)^2+\l(\f x2y\r)^2.$$
As $x^2=2m+1-y^2\eq 2m\ (\mo\ 8)$, $m$ is even and
$$\f m2\eq\l(\f x2\r)^2\eq\f x2\eq\f x2y\ \ (\mo\ 2).$$
We are done.
\qed

\heading{3. Proofs of Theorems 2 and 3}\endheading

\noindent{\tt Proof of Theorem 2}.
We distinguish four cases.

{\it Case} 1. $a=c=1$. Write $8=x_0^2+by_0^2+t_{z_0}$ with
$x_0,y_0,z_0\in\Z$, then $y_0\not=0$ and hence $8\gs b$. Since
$x^2+5y^2+t_z\not=13$, $x^2+6y^2+t_z\not=47$ and
$x^2+7y^2+t_z\not=20$, we must have $b\in\{1,2,3,4,8\}$.

{\it Case} 2. $a=1$ and $c=2$. Write $5=x_0^2+by_0^2+2t_{z_0}$
with $x_0,y_0,z_0\in\Z$. Then $y_0\not=0$ and hence $5\gs b$.
Observe that $x^2+3y^2+2t_z\not=8$ and $x^2+5y^2+2t_z\not=19$.
Therefore $b\in\{1,2,4\}$.

{\it Case} 3. $a=1$ and $c\gs3$. Since $2=x^2+by^2+ct_z$ for some
$x,y,z\in\Z$, we must have $b\ls 2$. If $b=1$, then there are
$x_0,y_0,z_0\in\Z$ such that $3=x_0^2+y_0^2+ct_{z_0}\gs c$ and
hence $c=3$. But $x^2+y^2+3t_z\not=6$, therefore $b=2$. For some
$x,y,z\in\Z$ we have $5=x^2+2y^2+ct_z\gs c$. Since
$x^2+2y^2+3t_z\not=23$ and $x^2+2y^2+5t_z\not=10$, $c$ must be
$4$.

{\it Case} 4. $a>1$. As $b\gs a\gs 2$ and $ax^2+by^2+ct_z=1$ for
some $x,y,z\in\Z$, we must have $c=1$. If $a>2$, then
$ax^2+by^2+t_z\not=2$. Thus $a=2$. For some $x_0,y_0,z_0\in\Z$ we
have $4=2x_0^2+by_0^2+t_{z_0}\gs b$. Note that
$2x^2+3y^2+t_z\not=7$ and $2x^2+4y^2+t_z\not=20$. Therefore $b=2$.

 In view of the above, Theorem 2 has been proven. \qed

\medskip
\noindent{\tt Proof of Theorem 3}. Let us first consider the case
$c>1$. Since $1=ax^2+bt_y+ct_z$ for some $x,y,z\in\Z$, we must
have $a=1$. Clearly $x^2+bt_y+ct_z\not=2$ if $c\gs 3$. So $c=2$.
For some $x_0,y_0,z_0\in\Z$ we have $5=x_0^2+bt_{y_0}+2t_{z_0}\gs
b$. It is easy to check that $x^2+3t_y+2t_z\not=8$. Therefore
$b\in\{2,4,5\}$.

Below we assume that $c=1$. If $a$ and $b$ are both greater than
$2$, then $ax^2+bt_y+t_z\not=2$. So $a\ls 2$ or $b\ls 2$.

{\it Case} 1. $a=1$. For some $x_0,y_0,z_0\in\Z$ we have
$8=x_0^2+bt_{y_0}+t_{z_0}\gs b$. Note that $x^2+5t_y+t_z\not=13$
and $x^2+7t_y+t_z\not=20$. So $b\in\{1,2,3,4,6,8\}$.

{\it Case} 2. $a=2$. For some $x_0,y_0,z_0\in\Z$ we have
$4=2x_0^2+bt_{y_0}+t_{z_0}\gs b$. Thus $b\in\{1,2,4\}$ since
$2x^2+3t_y+t_z\not=7$.

{\it Case} 3. $a>2$. In this case $b\ls 2$. If $b=1$, then for
some $x_0,y_0,z_0\in\Z$ we have $5=ax_0^2+t_{y_0}+t_{z_0}\gs a$,
and hence $a=4$ since $3x^2+t_y+t_z\not=8$ and
$5x^2+t_y+t_z\not=19$. If $b=2$, then for some $x,y,z\in\Z$ we
have $4=ax^2+2t_y+t_z\gs a$ and so $a\in\{3,4\}$.

The proof of Theorem 3 is now complete. \qed

\heading{4. Some conjectures and related discussion}\endheading

In this section we raise three related conjectures.

\proclaim{Conjecture 1} Any positive integer $n$ is a sum of a
square, an odd square and a triangular number. In other words,
each natural number can be written in the form $x^2+8t_y+t_z$ with
$x,y,z\in\Z$.
\endproclaim

We have verified Conjecture 1 for $n\ls 15,000$. By Theorem
1(iii), Conjecture 1 is valid when $n\not=t_4,t_8,t_{12},\ldots$.

\proclaim{Conjecture 2} Each $n\in\N$ can be written in any of the following forms with $x,y,z\in\Z$:
$$x^2+3y^2+t_z,\ x^2+3t_y+t_z,\ x^2+6t_y+t_z,\ 3x^2+2t_y+t_z,\ 4x^2+2t_y+t_z.$$
\endproclaim

\proclaim{Conjecture 3} Every $n\in\N$ can be written in the form
$x^2+2y^2+3t_z\ ($with $x,y,z\in\Z)$ except $n=23$, in the form
$x^2+5y^2+2t_z\ ($or the equivalent form $5x^2+t_y+t_z)$ except
$n=19$, in the form $x^2+6y^2+t_z$ except $n=47$, and in the form
$2x^2+4y^2+t_z$ except $n=20$.
\endproclaim

Both Conjectures 2 and 3 have been verified for $n\ls 10,000$.

The second statement in Conjecture 3 is related to an assertion of Ramanujan
confirmed by Dickson [D1] which states that
even natural numbers not of the form $4^k(16l+6)$ (with $k,l\in\N$)
can be written as $x^2+y^2+10z^2$ with $x,y,z\in\Z$.
Observe that
$$\align &n=x^2+5y^2+2t_z\ \t{for some}\ x,y,z\in\Z
\\\iff& 4n+1=x^2+5y^2+z^2\ \t{for some}\ x,y,z\in\Z\ \t{with}\ 2\nmid z
\\\iff&8n+2=2(x^2+y^2)+10z^2=(x+y)^2+(x-y)^2+10z^2\
\\&\qquad\qquad\qquad\qquad\qquad\ \t{for some}\ x,y,z\in\Z\ \t{with}\ 2\nmid y
\\\iff&8n+2=x^2+y^2+10z^2\ \t{for some}\ x,y,z\in\Z\ \t{with}\ x\not\eq y\ (\mo\ 4).
\endalign$$

Below we reduce the converses of Theorems 2 and 3 to Conjectures 1
and 2. For convenience, we call a ternary quadratic polynomial
$f(x,y,z)$ {\it essential} if $\{f(x,y,z):\,x,y,z\in\Z\}=\N$.
(Actually, in 1748 Goldbach (cf. [D2, p.\,11]) already stated that
$x^2+y^2+2t_z$, $x^2+2y^2+t_z$, $x^2+2y^2+2t_z$ and
$2x^2+2t_y+t_z$ are essential.)
\medskip

{\tt Step I}. We show that the 10 quadratic polynomials listed in
Theorem 2 are essential except for the form $x^2+3y^2+t_z$
appearing in Conjecture 2.

As $4x^2+y^2+2t_z\sim 4x^2+t_y+t_z$, the form $x^2+(2y)^2+2t_z$ is
essential by Theorem 1(i). Both $x^2+(2y)^2+t_z$ and
$2x^2+2y^2+t_z=(x+y)^2+(x-y)^2+t_z$ are essential by Theorem 1(ii)
and the trivial fact $t_z=0^2+0^2+t_z$. We have pointed out in
Section 1 that $x^2+2(2y)^2+t_z$ is essential by [JP, Theorem 6],
and we don't have an easy proof of this deep result.

 Since
$$x^2+2y^2+4t_z\sim x^2+2(t_y+t_z)\sim t_x+t_y+2t_z,$$
the form $x^2+2y^2+4t_z$ is essential by Theorem B (of Liouville).
By the Gauss-Legendre theorem, for each $n\in\N$ we can write
$8n+2=(4x)^2+(2y+1)^2+(2z+1)^2$ (i.e., $n=2x^2+t_y+t_z)$ with
$x,y,z\in\Z$. Thus the form $x^2+2y^2+2t_z$ is essential since
$2x^2+y^2+2t_z\sim 2x^2+t_y+t_z$.
\medskip

{\tt Step II}. We analyze the 15 quadratic polynomials listed in Theorem
3.

By Theorem 1(i), $(2x)^2+t_y+t_z$ and $x^2+t_y+t_z$ are essential.
Since
$$\gather x^2+2t_y+t_z\sim t_x+t_y+t_z,\\x^2+2t_y+2t_z\sim t_x+t_y+2t_z,
\\x^2+4t_y+2t_z\sim t_x+4t_y+t_z,
\\x^2+5t_y+2t_z\sim t_x+5t_y+t_z,
\\2x^2+4t_y+t_z\sim 2t_x+2t_y+t_z,
\endgather$$
the forms
$$x^2+2t_y+t_z,\ x^2+2t_y+2t_z,\ x^2+4t_y+2t_z,\ x^2+5t_y+2t_z,\ 2x^2+4t_y+t_z$$
are all essential by Liouville's theorem. For $n\in\N$ we can
write $2n=x^2+4t_y+2t_z$ with $x,y,z\in\Z$, and hence
$n=2x_0^2+2t_y+t_z$ with $x_0=x/2\in\Z$. So the form
$2x^2+2t_y+t_z$ is also essential.

Recall that $2x^2+t_y+t_z$ and $2x^2+y^2+2t_z$ are essential by
the last two sentences of Step I. For each $n\in\N$ we can choose
$x,y,z\in\Z$ such that $2n+1=2x^2+(2y+1)^2+2t_z$ and hence
$n=x^2+4t_y+t_z$. So the form $x^2+4t_y+t_z$ is essential.

The remaining forms listed in Theorem 3 are $x^2+8t_y+t_z$ and
four other forms, which appear in Conjectures 1 and 2
respectively. We are done.
\medskip

\Ack. The author is indebted to the referee for his/her helpful comments.

\medskip
\noindent {\bf Added in proof}. The second conjecture in Section 4
has been confirmed by Song Guo, Hao Pan and the author.

\bigskip

\leftline{Department of Mathematics} \leftline{Nanjing University}
\leftline{Nanjing 210093} \leftline{People's Republic of China}
\leftline {\tt zwsun\@nju.edu.cn}
\leftline {\tt http://math.nju.edu.cn/${}^\sim$zwsun}

\enddocument